\documentclass[12pt]{article}
\usepackage{amsmath,amssymb,amsthm}

\newtheorem{theorem}{Theorem}

\newtheorem{corollary}[theorem]{Corollary}

\newtheorem{conjecture}[theorem]{Conjecture}

\newcommand\R{\mathbb{R}}
\newcommand\Z{\mathbb{Z}}
\newcommand\Prb{\mathbb{P}}
\newcommand\bz{\mathbf{z}}
\newcommand\cA{\mathcal{A}}
\newcommand\cB{\mathcal{B}}
\newcommand\cC{\mathcal{C}}
\newcommand\cM{\mathcal{M}}
\newcommand\Th[1]{Theorem~\ref{t:#1}}
\newcommand\Co[1]{Corollary~\ref{c:#1}}
\newcommand\Fg[1]{Figure~\ref{f:#1}}
\newcommand\Sc[1]{Section~\ref{s:#1}}

\newdimen\unit\newdimen\psep\newcount\nd\newcount\ndx\newbox\dotb\newbox\ptbox
\newdimen\dx\newdimen\dy\newdimen\dxx\newdimen\dyy\newdimen\hgt
\newdimen\xoff\newdimen\yoff
\newcommand\clap[1]{\hbox to 0pt{\hss{#1}\hss}}
\newcommand\vdisk[1]{{\font\dotf=cmr10 scaled #1\dotf.}}
\newcommand\varline[2]{\setbox\dotb\hbox{\vdisk{#1}}\xoff=-.5\wd\dotb
\wd\dotb=0pt\yoff=-.5\ht\dotb\psep=#2\ht\dotb}
\newcommand\varpt[1]{\setbox\ptbox\clap{\vdisk{#1}}\setbox\ptbox
\hbox{\raise-.5\ht\ptbox\box\ptbox}}
\newcommand\cpt{\copy\ptbox}
\newcommand\point[3]{\rlap{\kern#1\unit\raise#2\unit\hbox{#3}}}
\newcommand\setnd[4]{\dx=#3\unit\advance\dx-#1\unit\divide\dx by\psep
\dy=#4\unit\advance\dy-#2\unit\divide\dy by\psep
\multiply\dx by\dx\multiply\dy by\dy\advance\dx\dy\nd=1\advance\dx-1sp
\loop\ifnum\dx>0\advance\dx-\nd sp\advance\nd1\advance\dx-\nd sp\repeat}
\newcommand\dl[4]{{\setnd{#1}{#2}{#3}{#4}\dline{#1}{#2}{#3}{#4}\nd}}
\newcommand\dline[5]{{\nd=#5\hgt=#2\unit\dx=#3\unit\advance\dx-#1\unit
\divide\dx by\nd\dy=#4\unit\advance\dy-#2\unit\divide\dy by\nd
\advance\hgt\yoff\rlap{\kern#1\unit\kern\xoff\loop\ifnum\nd>1\advance\nd-1
\advance\hgt\dy\kern\dx\raise\hgt\copy\dotb\repeat}}}

\newcommand\qellip[4]{{\setnd{0}{0}{#3}{#4}\dx=\unit\dy=0pt\raise\yoff\rlap{%
\kern#1\unit\kern\xoff\raise#2\unit\hbox{\loop\ifnum\dx>0\rlap{\kern#3\dx
\raise#4\dy\copy\dotb}\hgt=\dx\divide\hgt by\nd\advance\dy\hgt\hgt=\dy
\divide\hgt by\nd\advance\dx-\hgt\repeat\rlap{\raise#4\dy\copy\dotb}}}}}
\newcommand\bez[6]{{\setnd{#1}{#2}{#3}{#4}\ndx=\nd\setnd{#3}{#4}{#5}{#6}
\ifnum\ndx>\nd\nd=\ndx\fi\dx=#3\unit\advance\dx-#1\unit\dy=#4\unit
\advance\dy-#2\unit\dxx=#5\unit\advance\dxx-#1\unit\dyy=#6\unit\advance
\dyy-#2\unit\advance\dxx-2\dx\advance\dyy-2\dy\divide\dxx by\nd\divide\dyy
by\nd\advance\dx.25\dxx\advance\dy.25\dyy\divide\dx by\nd\divide\dy by\nd
\multiply\nd by2\dx=100\dx\dy=100\dy\dxx=100\dxx\dyy=100\dyy\divide\dxx by\nd
\divide\dyy by\nd\hgt=#2\unit\raise\yoff\rlap{\kern#1\unit\kern\xoff
\raise\hgt\copy\dotb\loop\ifnum\nd>0\advance\nd-1\advance\hgt0.01\dy
\kern0.01\dx\raise\hgt\copy\dotb\advance\dx\dxx\advance\dy\dyy\repeat}}}
\newcommand\ptu[3]{\point{#1}{#2}{\cpt\raise1ex\clap{$\scriptstyle{#3}$}}}
\newcommand\ptd[3]{\point{#1}{#2}{\cpt\raise-1.8ex\clap{$\scriptstyle{#3}$}}}
\newcommand\ptr[3]{\point{#1}{#2}{\cpt\raise-.4ex\rlap{$\ \scriptstyle{#3}$}}}
\newcommand\ptl[3]{\point{#1}{#2}{\cpt\raise-.4ex\llap{$\scriptstyle{#3}\ $}}}
\newcommand\ptlu[3]{\point{#1}{#2}{\raise.8ex\clap{$\scriptstyle{#3}$}}}
\newcommand\ptld[3]{\point{#1}{#2}{\raise-1.6ex\clap{$\scriptstyle{#3}$}}}
\newcommand\ptlr[3]{\point{#1}{#2}{\raise-.4ex\rlap{$\,\scriptstyle{#3}$}}}
\newcommand\ptll[3]{\point{#1}{#2}{\raise-.4ex\llap{$\scriptstyle{#3}\,$}}}

\newcommand\pnt[3]{\point{#1}{#2}{\clap{$#3$}}}

\newcommand\thnline{\varline{400}{.4}}

\varpt{2500}\thnline\unit=2em
\begin{document}
\title{Projections, Entropy and Sumsets}

\author{Paul Balister\thanks{Department of Mathematical Sciences,
University of Memphis, Memphis TN 38152, USA} \and B\'ela
Bollob\'as$^*$\thanks{Department of Pure Mathematics and
Mathematical Statistics, University of Cambridge, UK}
\thanks{Research supported in part by NSF grants CCR-0225610,
DMS-0505550 and W911NF-06-1-0076}}
\maketitle

\begin{abstract}
In this paper we have shall generalize Shearer's entropy inequality
and its recent extensions by Madiman and Tetali, and shall apply
projection inequalities to deduce extensions of some of the
inequalities concerning sums of sets of integers proved recently by
Gyarmati, Matolcsi and Ruzsa. We shall also discuss projection and
entropy inequalities and their connections.
\end{abstract}

\section{Introduction}

In 1949, Loomis and Whitney~\cite{LW} proved a fundamental inequality bounding
the volume of a body in terms of its $(n-1)$-dimensional projections. Over
forty years later, this inequality was extended considerably by Bollob\'as and
Thomason~\cite{BT}: they showed that a certain `box' is a solution of much
more general isoperimetric problems.

In 1978, Han~\cite{Han} proved the exact analogue of the
Loomis-Whitney inequality for the entropy of a family
$\{X_1,\dots,X_n\}$ of random variables, and in the same year
Shearer proved (implicitly) a considerable extension of this
inequality, namely the entropy analogue of the projection inequality
that was to be used some years later in \cite{BT} to deduce the Box
Theorem. (This extension was published only in 1986, in
\cite{CGFS}.) Recently, Madiman and Tetali~\cite{MT1, MT2}
strengthened Shearer's inequality to a two-sided inequality
concerning the joint entropy $H(X_1,\dots,X_n)$.

In this paper we have two main aims. The first is to prove an
entropy inequality that extends {\em both} sides of the
Madiman-Tetali inequality. Surprisingly, this inequality is not only
{\em much} more general than the earlier inequalities, but is also
just about {\em trivial}. Our second aim is to point out that the
projection inequalities imply extensions of some very recent
inequalities of Gyarmati, Matolcsi and Ruzsa~\cite{GyMR} concerning
sums of sets of integers.

Our paper is organized as follows. In the next two sections we shall
review some of the projection and entropy inequalities. In \Sc{new}
we shall prove our extremely simple but very general entropy
inequality extending those of Shearer, and Madiman and Tetali. In
\Sc{sum} we shall turn to sumsets, and continue the work of
Gyarmati, Matolcsi and Ruzsa. Finally, in \Sc{prob}, we shall state
some related unsolved problems.

\section{Projection inequalities}

As in \cite{BT}, we call a compact subset of $\R^n$ which is the closure of
its interior a {\em body}, and write $\{e_1,\dots,e_n\}$ for the canonical
basis of $\R^n$. Given a body $K\subseteq\R^n$ and a set
$A\subseteq [n]=\{1,\dots,n\}$ of $d$ indices, we denote by $K_A$ the
orthogonal projection of $K$ into the linear span of the vectors $e_i$,
$i\in A$, and write $|K_A|$ for its $d$-dimensional Euclidean volume. (In
particular, $K_{[n]}=K$.) The volumes $|K_A|$ can be viewed as a measure of
the `perimeter' of $K$. In 1949, Loomis and Whitney~\cite{LW} (see also
\cite{A}, \cite[page 95]{BZ} and \cite[page 162]{H}) proved the following
isoperimetric inequality:
\begin{equation}\label{LuWh}
 |K|^{n-1}\le\prod_{i=1}^n|K_{[n]\setminus \{i\}}|.
\end{equation}
Close to fifty years later, Bollob\'as and Thomason~\cite{BT} proved
the following {\em Box Theorem\/} showing that for the {\em set\/}
of projection volumes $|K_A|$, $A\subseteq [n]$, the solution of the
isoperimetric problem is a {\em box}, i.e., a rectangular parallelepiped
whose sides are parallel to the coordinate axes.

\begin{theorem}\label{t:BTBT}
 Given a body\/ $K\subseteq\R^n$, there is a box\/ $B\subseteq\R^n$ with\/
 $|K|=|B|$ and\/ $|K_A|\ge|B_A|$ for every\/ $A\subseteq[n]$.
 \hfill\qed
\end{theorem}
This theorem is equivalent to the assertion that there exist constants
$k_i\ge0$ such that
\begin{equation}
 |K|=\prod_{i=1}^nk_i\quad\text{and}\quad
 |K_A|\ge\prod_{i\in A}k_i\text{ for all }A\subseteq[n].
\end{equation}
An immediate consequence of \Th{BTBT} is that, if the volume of a box can be
bounded in terms of the volumes of a certain collection of projections, then
the same bound will be valid for all bodies. In particular, the Loomis-Whitney
Inequality \eqref{LuWh} is an immediate consequence of the Box Theorem. In
fact, the Box Theorem was deduced from the Uniform Cover Inequality, which is
an even more obvious extension of~\eqref{LuWh}. To state this inequality, we
call a multiset $\cA$ of subsets of $[n]$ such that each element $i\in [n]$ is
in at least $k$ of the members of $\cA$ a $k$-{\em cover} of $[n]$. A
$k$-{\em uniform cover} or {\em uniform} $k$-{\em cover} is one in which every
element is in precisely $k$ members of~$\cA$. Thus the sets
$[n]\setminus\{i\}$ appearing in the Loomis-Whitney inequality \eqref{LuWh}
form an $(n-1)$-uniform cover of $[n]$. The Uniform Cover Inequality states
that if $K$ is a body in $\R^n$ and $\cA$ is a $k$-uniform cover of $[n]$ then
\begin{equation}\label{UC}
 |K|^k \le \prod_{A\in\cA} |K_A|.
\end{equation}
Clearly, the Uniform Cover Inequality is a trivial consequence of the Box
Theorem. Uniformity {\em is} needed for \eqref{UC} to hold: if $\cA$ is not
$k$-uniform, then \eqref{UC} does not hold for every body $K$, not even if
$\cA$ is a $k$-cover. Indeed, if $|K_A|<1$ for some $A$, then we can add an
arbitrary number of copies of $A$ to $\cA$, making the right hand side of
\eqref{UC} arbitrarily small.

By identifying a lattice point $\bz\in\Z^n$ with the unit cube
$Q_{\bz}\subseteq\R^n$ with centre $\bz$, \eqref{UC} implies that if $S$ is a
finite subset of $\Z^n$ and $S_A$ is the projection of $S$ to the subspace
spanned by $\{e_i\colon i\in A\}$, then for every uniform $k$-cover $\cA$ of
$[n]$ we have
\begin{equation}\label{UC-latt}
 |S|^k\le \prod_{A\in\cA}|S_A|.
\end{equation}
In fact, in this inequality we do not have to demand that the $k$-cover
$\cA=\{A_i\}$ is uniform: if $A'\subseteq A$ then $|S_{A'}|\le|S_A|$;
therefore, by removing elements from the sets $A_i$ so as to obtain a
{\em uniform} $k$-cover $\cA'=\{A_i'\}$ with $A'_i\subseteq A_i$, we have
$|S|^k\le \prod_i|S_{A'_i}|\le\prod_i|S_{A_i}|$.

\section{Entropy Inequalities}

Let us turn to some entropy inequalities related to the projection
inequalities above. As usual, we write $H(X)$ for the entropy of a random
variable $X$; in particular, if $X$ is a discrete random variable, then
\begin{equation*}
 H(X)=-\sum_k \Prb(X=k)\log_2\Prb(X=k).
\end{equation*}
It is easily seen that if $X$ takes $n$ values then $H(X)\le\log_2 n$, with
equality if and only if $X$ is uniformly distributed, i.e., takes every value
with probability $1/n$. If $X$ and $Y$ are two discrete random variables, then
the entropy of $X$ conditional on $Y$ is
\begin{equation*}
 H(X\mid Y)=-\sum_{k,l}\Prb(X=k,Y=l)\log_2\Prb(X=k\mid Y=l).
\end{equation*}

The entropy satisfies the following basic inequalities:
\begin{gather}
\label{bayes} H(X,Y)=H(X\mid Y)+H(Y),\\
\label{cond}  0\le H(X\mid Y)\le H(X),\\
\label{mono}  H(X\mid Y,Z)\le H(X\mid Y),
\end{gather}
where, for example, we write $H(X,Y)$ for the entropy of the joint
variable $(X,Y)$.

Analogously to our notation concerning projections, given a sequence
$X=(X_1,\dots,X_n)$ of $n$ random variables, for $A\subseteq [n]$ we
write $X_A=(X_i)_{i\in A}$. In 1978 Shearer proved the following
analogue of \eqref{UC} for entropy (the result was first published
in~\cite{CGFS}). Since $H(X_A)$ is a monotone increasing function of
$A$, in this inequality it makes no difference whether we take $\cA$
to be a $k$-cover or uniform $k$-cover.
\begin{theorem}\label{t:Shear}
 If\/ $\cA$ is a uniform $k$-cover of\/ $[n]$ then
 \begin{equation}\label{SI}
  kH(X)\le \sum_{A\in\cA}H(X_A).
 \end{equation}
\end{theorem}
A little earlier Han~\cite{Han} had proved the `Loomis-Whitney' form of
\Th{Shear}: $(n-1)H(X)\le \sum_i H(X_{[n]\setminus\{i\}})$. The first
non-trivial case of this inequality is $2H(X,Y,Z)\le H(X,Y)+H(X,Z)+H(Y,Z)$.
Curiously, in \cite{CGFS} it is remarked that this special case can be proved
analogously to what we stated as \Th{Shear}, and so can the case when $\cA$ is
the collection of all $k$-subsets of $[n]$.

Some years after the publication of \cite{BT} it was noted that \Th{Shear}
implies \Th{BTBT}. In fact, the reverse implication is also easy: this follows
from the fact that if $p_1,\dots,p_n$ are fixed `probabilities' with
$\sum p_i=1$ and $Np_i$ is an integer for every $i$, then the number of
sequences of length $N$ with $Np_i$ terms equal to $i$ is $2^{(1+o(1))H(X)N}$,
where $X$ is a random variable with $\Prb(X=i)=p_i$. Given random variables
$X_1,\dots,X_n$, we may assume that $X_i$ takes values in $V_i\subseteq \Z$,
so that $X=(X_1,\dots,X_n)$ takes values in $V=V_1\times\dots\times V_n$, and
there is an integer $d$ such that $d\,\Prb(X=v)$ is an integer for every
$v\in V$. Let $N$ be a multiple of $d$, and let
$S\subseteq V^N\subseteq\Z^{nN}$ be the set of all sequences in which $v$
occurs precisely $N\,\Prb(X=v)$ times. For $A\subseteq[n]$, write
$\tilde A\subseteq[nN]$ for the set of all coordinates of
$V^N\subseteq\Z^{nN}$ that correspond to one of the factors $V_i$, $i\in A$.
Then $S_{\tilde A}$ is the set of sequences in $V_A^N$ where each value
$v\in V_A$ occurs $N\,\Prb(X_A=v)$ times. If $\cA$ is a $k$ uniform cover of
$[n]$ then $\tilde\cA=\{\tilde A:A\in\cA\}$ is a $k$ uniform cover of $[nN]$
and so by \Th{BTBT}
\begin{equation*}
 |S|^k\le\prod_{A\in\cA}|S_{\tilde A}|.
\end{equation*}
Thus
\begin{equation*}
 2^{k(1+o(1))H(X)N}\le \prod_{A\in\cA}2^{(1+o(1))H(X_A)N}
\end{equation*}
and \Th{Shear} follows by letting $N\to\infty$.

Recently, Madiman and Tetali~\cite{MT1}, \cite{MT2} strengthened \Th{Shear} by
replacing the entropies $H(X_A)$ by certain conditional entropies; furthermore,
they also gave lower bounds for $H(X)$.

\begin{theorem}\label{t:MT}
 Let\/ $X=(X_i)_1^n$ be a sequence of random variables with\/ $H(X)$
 finite, and\/ $\cA$ a uniform\/ $k$-cover of\/ $[n]$. For\/
 $A\subseteq[n]$ with minimal element\/ $a\ge 1$ and maximal element\/ $b$,
 set\/ $A_*=\{1,\dots,a-1\}$ and\/ $A^*=\{i\notin A: 1\le i\le b-1\}$. Then
 \[
  \sum_{A\in\cA} H(X_A\mid X_{A^{^*}})\le k H(X)
  \le \sum_{A\in\cA} H(X_A\mid X_{A_*}).
  \eqno\qed
 \]
\end{theorem}
It should be noted that Theorem~\ref{t:MT} does {\em not\/} follow
from Shearer's Inequality, Theorem~\ref{t:Shear}.

Trivially, in the lower bound $\cA$ may be replaced by a
$k$-packing or a fractional $k$-packing, and in the upper bound it
may be replaced by a $k$-cover or a fractional $k$-cover, with the
obvious definitions.

\section{New Entropy Inequalities}\label{s:new}

Since, as shown in \cite{BT}, the Box Theorem follows from the
Uniform Cover Inequality \eqref{UC}, one has a Box Theorem type
strengthening of Shearer's Inequality; in fact, there is a similar
strengthening of \Th{MT} as well.

\begin{theorem}\label{t:MT-BI}
 Let\/ $X=(X_i)_1^n$ be a sequence of random variables with $H(X)$
 finite. Then there are non-negative constants $h_1,\dots,h_n$
 such that $H(X)=\!\sum_i^n h_i$ and
 \begin{equation*}
  H(X_A\mid X_{A^*})\le\sum_{i\in A}h_i
  \le H(X_A\mid X_{A_*})\quad\text{for all}\quad A\subseteq[n].
 \end{equation*}
\end{theorem}
\begin{proof}
We may take $h_i=H(X_i\mid X_{[i-1]})$; to prove the inequalities,
we inductively apply properties (\ref{bayes}--\ref{mono}).
\end{proof}

Although \Th{MT} does not follow from \Th{Shear} (Shearer's
Inequality), as we shall see now, it does follow from a result which
is extremely easy to prove but is still a considerable extension of
Shearer's Inequality and a generalization of the submodularity of
the entropy. Before we state this new inequality, we shall recall a
consequence of the basic entropy inequalities, and introduce a
partial order on the collection of multisets of subsets of $[n]$.

First, from \eqref{mono} and \eqref{bayes} one can deduce that
$H(X_A)$ is a {\em submodular} function of the set $A$: if
$A,B\subseteq[n]$ then
\begin{equation}\label{submod}
 H(X_{A\cup B})+H(X_{A\cap B})\le H(X_A)+H(X_B).
\end{equation}
To see this, note that by \eqref{mono} we have
\begin{equation*}
 H(X_{B\setminus A}\mid X_A)\le H(X_{B\setminus A}\mid X_{A\cap B});
\end{equation*}
using \eqref{bayes} to expand the first and last terms, we get
\begin{equation*}
 H(X_{A\cup B})-H(X_A)\le H(X_B)-H(X_{A\cap B}),
\end{equation*}
which is \eqref{submod}.

\bigskip
Second, let $\cM_{n,m}$ be the family of multisets of non-empty subsets of
$[n]$ with a total of $m$ elements. Given a multiset
$\cA=\{A_1,\dots,A_{\ell}\}\in\cM_{n,m}$ with non-nested sets $A_i$ and
$A_j$ (thus neither $A_i\subseteq A_j$ nor $A_j\subseteq A_i$ holds),
let $\cA'=\cA_{(ij)}$ be obtained from $\cA$ by replacing $A_i$ and $A_j$ by
$A_i\cap A_j$ and $A_i\cup A_j$, keeping only $A_i\cup A_j$ if
$A_i\cap A_j=\emptyset$. (If $A_i$ and $A_j$ are nested then replacing
$(A_i, A_j)$ by $(A_i\cap A_j,A_i\cup A_j)$ does not change $\cA$.) We call
$\cA'$ an {\em elementary compression} of $\cA$. Also, we call the
result of a sequence of elementary compressions a {\em compression}.

Let us define a partial order on $\cM_{(n,m)}$ by setting $\cA>\cB$ if $\cB$
is a compression of $\cA$. That $`>'$ defines a partial order on $\cM_{n,m}$
follows from the fact that if $\cA'$ is an elementary compression of $\cA$
then
\begin{equation*}
 \sum_{A\in\cA}|A|^2<\sum_{A\in\cA'}|A|^2.
\end{equation*}
Note that for every multiset $\cA\in\cM_{(n,m)}$ there is a unique minimal
multiset $\cA^\sharp$ dominated by $\cA$ consisting of the sets
\begin{equation*}
 A^\sharp_j=\{i\in[n]\mid\, i
 \text{ lies in at least $j$ of the sets }A\in\cA\}.
\end{equation*}
Equivalently, $\cA^\sharp$ is the unique multiset that is totally ordered
by inclusion and has the same multiset union as $\cA$.
If we renumber $[n]$ in such a way that each $A^\sharp_j$ is an initial
segment, then the $A^\sharp_j$s are just the rows of the Young tableaux
associated with the set system $\cA$, as in \Fg{YT}. In particular, if
$\cA$ is a $k$-uniform cover, then $\cA^\sharp=k\{[n]\}$.

\begin{figure}
\[\unit=20pt
 \dl0005\dl1014\dl2024\dl3033\dl4041\dl5051\dl6061
 \dl0070\dl0161\dl0232\dl0333\dl0424
 \point{-.5}{0.3}{$
  \pnt20{\circ}   \pnt30{\circ}   \pnt40{\circ}
  \pnt10{\bullet} \pnt31{\bullet} \pnt50{\bullet}
  \pnt11{\diamond}\pnt21{\diamond}\pnt60{\diamond}
  \pnt12{\ast}    \pnt22{\ast}
  \pnt13{+}       \pnt32{+}       \pnt23{\times}
  \pnt{1}{-1}{1}  \pnt{2}{-1}{2}  \pnt{3}{-1}{3}
  \pnt{4}{-1}{4}  \pnt{5}{-1}{5}  \pnt{6}{-1}{6}
  $}
 \point{8}{4.4}{$\cA^\sharp=$}
 \point{8}{3.4}{$\{1,2\}$}
 \point{8}{2.4}{$\{1,2,3\}$}
 \point{8}{1.4}{$\{1,2,3\}$}
 \point{8}{0.4}{$\{1,2,3,4,5,6\}$}
 \hskip11\unit
\]
\caption{The minimal compression $\cA^\sharp$ of
 $\cA=\{\{2,3,4\},\{1,3,5\},\{1,2,6\},$ $\{1,2\},\{1,3\},\{2\}\}$.}
 \label{f:YT}
\end{figure}

Here is then our essentially trivial but general entropy inequality.
\begin{theorem}\label{t:gen1}
 Let\/ $X=(X_i)_1^n$ be a sequence of random variables with\/ $H(X)$
 finite, and let\/ $\cA$ and\/ $\cB$ be finite multisets of subsets
 of\/ $[n]$. If\/ $\cA>\cB$ then
 \begin{equation}\label{gen-ineq}
  \sum_{A\in\cA}H(X_A)\ge\sum_{B\in\cB}H(X_B).
 \end{equation}
\end{theorem}
\begin{proof}
All we have to check is that \eqref{gen-ineq} holds if $\cB$ is an
elementary compression of $\cA$, i.e., if $\cB=\cA_{(ij)}$ for some
$i$ and $j$, where $\cA=\{A_1,\dots,A_{\ell}\}$. But then \eqref{gen-ineq}
is equivalent to
\begin{equation*}
 H(X_{A_i})+H(X_{A_j})\ge H(X_{A_i\cap A_j})+H(X_{A_i\cup A_j}),
\end{equation*}
which holds by \eqref{submod}, the submodularity of the entropy.
\end{proof}

We dignify the special case of Theorem~\ref{t:gen1} in which
${\mathcal B}$ is the minimal multiset ${\mathcal A}^{\sharp}$
dominated by ${\mathcal A}$ by calling it a theorem. This is the
inequality one is most likely to use.

\begin{theorem}\label{t:gen2}
 Let\/ $X=(X_i)_1^n$ be a sequence of random variables with\/ $H(X)$
 finite, and let\/ $\cA$ be a finite multiset of subsets of\/ $[n]$. Then
 \[
  \sum_{A\in\cA^\sharp}H(X_A)\le\sum_{A\in\cA}H(X_A).
  \eqno\qed
 \]
\end{theorem}

Let us illustrate \Th{gen1} with a simple example: as
$\{\{1,2\},\{1,3\},\{4\}\}>\{\{1,2,3\},\{1,4\}\}$,
\begin{equation*}
 H(X_1,X_2)+H(X_1,X_3)+H(X_4)\ge H(X_1,X_2,X_3)+H(X_1,X_4).
\end{equation*}
Also, let us point out that even Theorem~\ref{t:gen2} is stronger
than \Th{MT}, the Madiman--Tetali inequality.
\begin{proof}[Proof of Theorem~\ref{t:gen2} $\Rightarrow$ \Th{MT}.]
Since $H(X_A\mid X_B)=H(X_{A\cup B})-H(X_B)$, the upper bound inequality is
\begin{equation*}
 k H(X)+\sum_{A\in\cA}H(X_{A_*})\le \sum_{A\in\cA} H(X_{A\cup A_*}),
\end{equation*}
which follows from the fact that the multiset
$\cC_1=\{A_*:A\in\cA\}\cup k\{[n]\}$ is totally ordered and has the same
multiset union as $\cC_2=\{A\cup A_*:A\in\cA\}$, so $\cC_1=\cC_2^\sharp$.
Similarly, the lower bound inequality is equivalent to
\begin{equation*}
 \sum_{A\in\cA}H(X_{A\cup A^*})\le \sum_{A\in\cA} H(X_{A^*})+kH(X).
\end{equation*}
which follows from the fact that the multiset
$\cC_3=\{A\cup A^*:A\in\cA\}$ is totally ordered and has the same
multiset union as $\cC_4=\{A_*:A\in\cA\}\cup k\{[n]\}$, so
$\cC_3=\cC_4^\sharp$.
\end{proof}

The inequality corresponding to Theorem~\ref{t:gen2} in terms of
projections of bodies is false. For example, consider the set $K$ in
\Fg{1}. Then $|K|=5$, $|K_{\{1\}}|=2$, but
$|K_{\{1,2\}}|=|K_{\{1,3\}}|=3$, so
\begin{equation*}
 |K_{\{1,2,3\}}||K_{\{1\}}|>|K_{\{1,2\}}||K_{\{1,3\}}|.
\end{equation*}

\begin{figure}
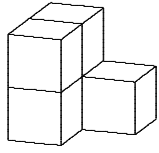

\begin{equation*}\unit=20pt
 \dl{0}{0.1}{1}{0}\dl{0}{1.1}{1}{1}\dl{0}{2.1}{1}{2}
 \dl{1.4}{0.3}{2.4}{0.2}\dl{1.4}{1.3}{2.4}{1.2}
 \dl{1.8}{1.6}{2.8}{1.5}\dl{0.4}{2.4}{1.4}{2.3}
 \dl{0.8}{2.7}{1.8}{2.6}
 \dl{0}{0.1}{0}{2.1}\dl{1}{0}{1}{2}
 \dl{1.4}{0.3}{1.4}{2.3}\dl{1.8}{1.6}{1.8}{2.6}
 \dl{2.4}{0.2}{2.4}{1.2}\dl{2.8}{0.5}{2.8}{1.5}
 \dl{1}{0}{1.4}{0.3}\dl{2.4}{0.2}{2.8}{0.5}
 \dl{1}{1}{1.8}{1.6}\dl{2.4}{1.2}{2.8}{1.5}
 \dl{1}{2}{1.8}{2.6}\dl{0}{2.1}{0.8}{2.7}
 \hskip3\unit
\end{equation*}
\caption{A body $K$ made up of five unit cubes. Coordinate `1' is
horizontal.}\label{f:1}
\end{figure}

\section{Sumsets}\label{s:sum}

Let $S_1,\dots,S_n$ be finite sets in a commutative semigroup with sum
\begin{equation*}
 S=S_1+\dots+S_n=\{s_1+\dots+s_n: s_i\in S_i\text{ for every }i\}.
\end{equation*}
For $A\subseteq[n]$ set $S_A=\sum_{i\in A}S_i$, so that $S_{[n]}=S$. We shall
think of $S$ as an $n$-dimensional body in $\R^n$ and $S_A$ as its canonical
projection into the subspace spanned by $\{e_i:i\in A\}$. Gyarmati, Matolcsi
and Ruzsa~\cite{GyMR} proved the analogue of the Loomis-Whitney inequality in
this context. In fact, the analogue of the Uniform Cover inequality and Box
Theorem are just as easy to show.

To see this, put an arbitrary linear order on each of the sets
$S_i$. For each $A=\{i_1,\dots,i_r\}\subseteq[n]$ define an
embedding $\varphi_A$ of $S_A$ into the Cartesian product
$\prod_{i\in A}S_i$ by mapping $s\in S_A$ to the lexicographically
least element $(s_{i_1}, \dots , s_{i_r})$ of $\prod_{i\in A}S_i$
with coordinates summing to $s$. (In fact, there are many other
orders we could choose instead of the lexicographic order: all we
need is that the assertions below hold for these orders.) As shown
by Gyarmati, Matolcsi and Ruzsa~\cite{GyMR}, the projection of
$S'=\varphi_{[n]}(S_{[n]})$ into $\prod_{i\in A}S_i$ is contained in
$\varphi_A(S_A)$. To see this, note that if $(s_1,\dots,s_n)\in S'$
then any projection $(s_{i_1},\dots,s_{i_r})$ is lexicographically
minimal with the same sum, since if
$s_{i_1}+\dots+s_{i_r}=s'_{i_1}+\dots+s'_{i_r}$ with
$(s'_{i_1},\dots,s'_{i_r})<(s_{i_1},\dots,s_{i_r})$, then
$s_1+\dots+s_n=s'_1+\dots+s'_n$ and
$(s'_1,\dots,s'_n)<(s_1,\dots,s_n)$ where $s'_i=s_i$ if $i\notin A$.
Thus $|(S')_A|\le|\varphi_A(S_A)|=|S_A|$. Now the following result
is immediate from \Th{BTBT} applied to~$S'$.

\begin{theorem}
  There are constants\/ $\lambda_1,\dots,\lambda_n\ge0$ such that
 \[
  |S|=\prod_1^n \lambda_i\quad\text{and}\quad
  |S_A|\ge \prod_{i\in A}\lambda_i \ \text{ for all }\ A\subseteq[n].
 \]
 In particular, if  $\cA$ is a uniform\/ $k$-cover\/ of\/ $[n]$ then
 \[
  |S|^k\le \prod_{A\in\cA}|S_A|.
  \eqno\qed
 \]
\end{theorem}

Using a similar approach, one can prove the following, which is
stated (with a slight error) as an open problem in~\cite{GyMR}.

\begin{theorem}\label{t:GyMR}
 If\/ $A,B_1,\dots,B_k$ are finite sets of integers and\/
 $C\subseteq B_1+\dots+B_k$, then
 \begin{equation}\label{e:cor}
  |A+C|^k\le |C|^{k-1}\prod_{i=1}^k|A+B_i|.
 \end{equation}
\end{theorem}
\begin{proof}
For convenience, write $n=k+1$ and $B_n=B_{k+1}=A$. Define maps
$\varphi_{T}$, $T\subseteq[n]$, as above. Let $S'$ be
$\varphi_{[n]}(C+B_n)$. Then $|S'_{[k]}|\le|C|$ and
$|S'_{\{i,n\}}|\le|B_i+B_n|$. The result follows by applying
\eqref{UC-latt} to $S'$ and the $k$-uniform cover of $[n]$
consisting of the following $2k-1$ sets: the $k$ pairs $\{1,n\}$,
$\{2,n\}$, $\dots$, $\{k,n\}$, and the $k$-set $[n-1]$ taken $k-1$
times.
\end{proof}

One can have equality in \Th{GyMR}, for example, if $A=[n]$,
$C=B_1=\{0,n\}$, $B_2=\dots=B_k=\{0\}$. This shows that inequality
\eqref{e:cor} may break down if $|C|^{k-1}$ is replaced by $|C|^i$
with $i<k-1$.

It is worth noting that one can prove a lower bound on $|S|$ which
is additive in the $|S_A|$ in the case when the sets $S_i$ lie in a
torsion free abelian group. This generalizes Theorem~1.1 of
\cite{GyMR}.

\begin{theorem}\label{t:add}
 If the sets\/ $S_i$ lie in a torsion-free abelian group then there
 are subsets\/ $S'_i\subseteq S_i$ of cardinality at most\/~$2$
 such that for any uniform $k$-cover $\cA$ of\/ $[n]$ we have
 \begin{equation*}
  k(|S|-1)\ge k(|S'|-1)\ge \sum_{A\in\cA} (|S_A|-1),
 \end{equation*}
 where $S'$ is the set of sums $s_1+\dots+s_k\in S$
 such that\/ $\{i:s_i\notin S'_i\}\subseteq A$ for some $A\in\cA$.
\end{theorem}
\begin{proof}
We first note that any torsion-free abelian group
can be given an ordering compatible with addition.

Pack a $k\times n$ grid with the sets $A\in\cA$ in the obvious
manner: each $A=\{j_1,\dots,j_r\}$ is packed as a set of pairs
$A'=\{(i_1,j_1),\dots,(i_r,j_r)\}$ so that the $A'$, $A\in\cA$, are
disjoint and cover the whole of $[k]\times[n]$. The $i_k$s are
otherwise arbitrarily chosen.

We may assume without loss of generality that the minimum
elements of $S_i$ are all equal to 0. Let $a_i$ be the maximum
element of $S_i$. The set $S'_i$ will be chosen to be
$\{0,a_i\}$. For convenience write $a_T=\sum_{i\in T}a_i$.
We shall mark $k$ copies of $S-\{0\}$ as follows.

Process each element of $[k]\times[n]$ in the lexicographic order
--- i.e.,
\begin{equation*}
 (1,1),\dots,(1,n),(2,1),\dots,(2,n),\dots\dots,(k,n).
\end{equation*}
Suppose we are processing $(i,j)$. Then $(i,j)=(i_t,j_t)$ for some
$A\in\cA$.
In the $i$'th copy of $S'-\{0\}$, mark all the elements that are in
\begin{equation*}
 a_{[j]-A}+S_A\cap(a_{[j-1]},a_{[j]}].
\end{equation*}
Note that all elements of $a_{[j]-A}+S_A$ lie in $S'$ (indeed in
$S_{A\cup[j]}\cap S'$), and, subtracting $a_{[j]-A}$, the number of elements
marked is equal to the number of elements of $S_A$ that lie in the
interval
\begin{equation*}
 (a_{[j-1]\cap A}, a_{[j]\cap A}].
\end{equation*}
(Note by assumption $j\in A$ so $[j]-A=[j-1]-A$).
Now it is clear that for distinct $(i,j)$, distinct elements are marked
(since they all lie in the $i$'th copy of $S'\cap (a_{[j-1]},a_{[j]}]$
and these sets are distinct), so at most $k(|S'|-1)$ elements are marked
in total. (The element $0\in S$ is not included in any of the intervals
$(a_{[j-1]},a_{[j]}]$.) However, every element in $S_A-\{0\}$ lies in
some interval $(a_{[j-1]\cap A},a_{[j]\cap A}]$ for some $j\in A$, so
results in some element being marked. Since it is clear that $|S|\ge|S'|$,
the result follows.
\end{proof}

\begin{corollary}\label{c:tfg}
 If the sets  $S_i$ lie in a torsion-free abelian group then there
 exists constants\/ $\sigma_i$ such that
 \[
  |S|-1=\sum_{i=1}^n\sigma_i\quad\text{and}\quad
  |S_A|-1\le\sum_{i\in A}\sigma_i\text{ for all }A\subseteq[n].
  \eqno\qed
 \]
\end{corollary}

\Th{add} fails for groups with torsion when, for example, all $S_i$
are equal to some non-trivial finite subgroup. If we insist that
$|S|$ is smaller than the order of the smallest non-trivial subgroup
then we have the famous Cauchy--Davenport theorem, which can be
written in the following form.

\begin{theorem}\label{t:cd}
 If\/ $S_1,\dots,S_n$ are non-empty subsets of\/ $\Z_p$ and\/
 $S=S_1+\dots+S_n$, then either\/ $|S|\ge p$ or
 \[
  \hspace{130pt} |S|-1\ge\sum_i(|S_i|-1). \hspace{130pt}{\square}
 \]
\end{theorem}

\Th{cd} is the analogue of \Co{tfg} for the 1-uniform cover
$\cA=\{\{1\},\dots,\{n\}\}$, and can be
extended to all finite (even non-abelian) groups as
is shown in \cite{K} and \cite{W} (see also~\cite{BW}).
\begin{theorem}\label{t:gcd}
 If\/ $S_1,\dots,S_n$ are non-empty subsets of a finite group $G$ and\/
 $S=S_1\star\dots\star S_n$ ($\star$ denoting the group operation),
 then either $|S|\ge p$ or
 \[
  |S|-1\ge\sum_i(|S_i|-1).
 \]
 where $p$ is the smallest prime dividing $|G|$. \hfill{$\square$}
\end{theorem}
Unfortunately, \Th{gcd} does not generalize to more general covers.
For example, if $S_1=S_2=S_3=\{0,1,3,5\}\subseteq \Z_{13}$ then
$|S_1+S_2|=|S_1+S_3|=|S_2+S_3|=9$ and $|S_1+S_2+S_3|=12$, so
\[
 2(|S_1+S_2+S_3|-1)<(|S_1+S_2|-1)+(|S_1+S_3|-1)+(|S_2+S_3|-1).
\]

\section{Conjectures}\label{s:prob}

The most obvious problems related to the results above concern
general (not necessarily commutative) groups. In fact, Ruzsa has
already asked whether a suitable analogue of the inequality
corresponding to the Loomis--Whitney inequality holds for all
groups. It is not unreasonable to hope that the analogue of the Box
Theorem (or Cover Inequality) holds as well, as does the extension
of Corollary~\ref{c:tfg}. To state these conjectures, given finite
non-empty sets $S_1, \dots , S_n$ in a group $G$ with operation
$\star$ as above, and a set $A\subset [n]$, write $N_A$ for the
maximal number of elements in a product set obtained from
$S_1\star\dots\star S_n$ by replacing each $S_i$, $i\notin A$, by a
single element of $S_i$. Similarly, write $n_A$ for the
corresponding minimum.
\begin{conjecture}
Let\/ $S_1,\dots,S_n$ be non-empty finite subsets of a group. Set\/
$S=S_1\star\dots\star S_n$, and let $N_A$ be as above. Then there
are constants\/ $\lambda_1,\dots,\lambda_n>0$ such that
 \[
  |S|=\prod_1^n \lambda_i\quad\text{and}\quad
  N_A\ge \prod_{i\in A}\lambda_i\ \ \text{for all}\ \ A\subseteq[n].
  \eqno\square
 \]
\end{conjecture}

\begin{conjecture}
Let\/ $S_1,\dots,S_n$ be non-empty finite subsets of a group, and
let $S$ and $n_A$ be as above. Then there are constants\/ $\sigma_i$
such that
 \[
  |S|-1=\sum_{i=1}^n\sigma_i\quad\text{and}\quad
  n_A-1\le\sum_{i\in A}\sigma_i\text{ for all }A\subseteq[n].
 \]
\end{conjecture}

In conclusion, we should say that both these conjectures are rather
tentative: we would not be amazed if they turned out to be false.

\section{Acknowledgements}

The results in \Sc{sum} were proved after (and while) listening to
I.~Ruzsa's lecture in Tel Aviv in June, 2007; we are grateful to
Professor Ruzsa for showing us his slides of this lecture, and for
prepublication access to~\cite{GyMR}.

\end{document}